\documentclass[12pt,a4paper]{article}

\usepackage[utf8]{inputenc}
\usepackage[russian]{babel}
\usepackage{amsthm}
\usepackage{amsfonts}
\usepackage{amssymb}
\usepackage{amsmath}
\usepackage[affil-it]{authblk}
\usepackage{mathrsfs}
\usepackage[usenames]{color}
\usepackage{graphicx}
\usepackage{cite}

\newtheorem{Theorem}{Теорема}

\newtheorem{Lemma}{Лемма}

\newenvironment{Def}{\par\bigskip\noindent{\bf Определение. }}{\par\bigskip}

\DeclareMathOperator{\re}{Re}

\DeclareMathOperator{\const}{const}

\newcommand{\CC}{\mathbb{C}}
\newcommand{\RR}{\mathbb{R}}
\newcommand{\NN}{\mathbb{N}}

\newcommand{\vint}[2]{\bigl|_{#1}^{#2}\bigr.}
\newcommand{\sint}[3]{\{#1\}_{#2}^{#3}}
\newcommand{\set}[2]{\bigl\{#1\,\bigl|\bigr.\,#2\bigr\}}

\makeatletter
\renewcommand{\section}{\@startsection{section}{1}%
{0pt}{3.5ex plus 1ex minus .2ex}%
{2.3ex plus.2ex}{\normalfont\large\bfseries}}

\makeatother

\title{О критерии Молчанова компактности резольвенты для несамосопряжённого оператора Штурма--Лиувилля}
\author{Туманов С.\,Н.}
\affil{Московский Центр фундаментальной и прикладной математики
при МГУ им. М.\,В. Ломоносова,
Москва, Россия}

\begin{document}

\maketitle

\begin{abstract}
Условие Молчанова является необходимым условием компактности резольвенты для широкого класса обыкновенных дифференциальных операторов
произвольного порядка, но уже для оператора Штурма--Лиувилля оно не является достаточным, даже
если вещественная часть потенциала неотрицательна. Критерий Молчанова остаётся справедливым в формулировке наиболее близкой к исходной для
потенциалов, принимающих значения в более узком секторе, чем полуплоскость, отделённом от отрицательной полуоси. Этим вопросам посвящена
настоящая работа.
\end{abstract}
\frenchspacing

\section{Введение}

Одним из существенных результатов теории качественного спектрального анализа обыкновенных дифференциальных операторов является Теорема Молчанова ---
критерий полной непрерывности резольвенты \cite{Molchanov}
сингулярного оператора Штурма--Лиувилля
$$
Ly=-y''+qy.
$$

Изначально сформулированный для операторов с вещественным потенциалом $q(x)\ge\const$ ($x\in\RR$ или $x\in\RR_+$), критерий получил обобщения
\cite{Brinck, Lidsky, Ismagilov, Birman, LevitanSuvorch, Fortunato, IsmagilovKost}, направленные, в основном, на ослабление условия
ограниченности $q$ снизу, либо на векторный случай, когда в роли $q$ выступает эрмитова матрица, $q\ge0$. Мы не затрагиваем
дифференциальных операторов в частных производных.

Обобщения, касающиеся комплекснозначного $q$, не многочисленны, и здесь ключевой является работа Лидского \cite{Lidsky},
где критерий Молчанова был обобщён для потенциалов, принимающих значения в одном из секторов комплексной плоскости:
$0\le\arg q(x)\le\pi/2$, либо $-\pi/2\le\arg q(x)\le0$.

Всюду далее
\begin{equation}
\label{eqqdefinit}
q\in L_{1,loc}(\RR_+)\ \mbox{--- комплекснозначная функция}.
\end{equation}

Будем говорить, что $q$ удовлетворяет {\it условию Молчанова}, если для любой бесконечной
системы непересекающихся интервалов $D_n\subset\RR_+$ равной длины
$$
\int\limits_{D_n}|q(x)|\,dx\to+\infty\ \mbox{при}\ n\to\infty.
$$

Уже для вещественного $q$ при отказе от ограниченности снизу, условие Молчанова перестаёт быть критерием полной непрерывности резольвенты. Например, для
$q(x)=-x^2$, см. \cite{Titchmarsh}.

Тем не менее, оно остаётся необходимым, даже если не накладывать на $q$ никаких ограничений сверх \eqref{eqqdefinit}.
Этот результат мы докажем
в том числе и для обыкновенных дифференциальных операторов произвольного порядка $n\ge 2$.

Дополняя его теоремами о достаточности и контрпримером, мы построим картину возможных обобщений Теоремы Молчанова на
комплекснозначный случай в терминах непосредственно условия Молчанова.

Рассмотрим дифференциальное выражение
\begin{equation}
\label{eqdiffexpr}
l(y)=-y''+qy
\end{equation}
и линеалы
\begin{gather*}
\mathscr{D}=\set{y\in L_2(\RR_+)}{y,y'\in AC_{loc}(\RR_+),\ l(y)\in L_2(\RR_+)},\\
\mathscr{D}_0=\set{y\in\mathscr{D}}{y(0)=y'(0)=0,\ \exists x_0>0\ \forall x\ge x_0\ y(x)=0},\\
\mathscr{D}_U=\set{y\in\mathscr{D}}{U(y)=0},
\end{gather*}
где $U$ --- некоторая форма краевых условий в $x=0$:
$$
U(y)=A y(0)+B y'(0),\quad A,B\in\CC,\ |A|+|B|>0.
$$

Определим дифференциальные операторы $L_0\subset L_U$ в $L_2(\RR_+)$ на соответствующих областях определения
$D_0\subset D_U$ дифференциальным выражением \eqref{eqdiffexpr}, следуя \cite{Naimark}.

\begin{Theorem}
\label{thnecessn2}
Для того, чтобы у оператора $L_0$ существовало расширение с компактной резольвентой, необходимо, чтобы $q$ удовлетворяло условию Молчанова.
\end{Theorem}

Этот результат является следствием общей теоремы \ref{thnecessity}, к которой мы обратимся в разделе \ref{sectnecessity}.

\begin{Def}
Скажем, что для потенциала $q$ выполнено $\RR^-$--{\it условие}, если при всех достаточно больших $x>x_0\ge0$ значения $q(x)$ лежат в секторе
$\alpha\le\arg (q(x)-q_0)\le\beta$ для некоторых $-\pi<\alpha\le\beta<\pi$ и $q_0\in\CC$.
\end{Def}

Другими словами, $\RR^-$--{\it условие} означает, что найдётся малый сектор, содержащий $\RR^-$, в котором $q-q_0$ асимптотически не принимает значений
для некоторого $q_0\in\CC$.

\begin{Def}
Потенциал $q$ назовём {\it секториальным}, если для него выполнено $\RR^-$--условие с $\beta-\alpha<\pi$.
\end{Def}

\begin{Theorem}
\label{thcriterion}
Пусть потенциал $q\in L_{1,loc}(\RR_+)$ секториальный. Тогда оператор $L_U$ имеет компактную резольвенту тогда и только тогда, когда $q$
удовлетворяет условию Молчанова.
\end{Theorem}

Как показывает следующая теорема, условие $\beta-\alpha<\pi$ не может быть ослаблено.
\begin{Theorem}
\label{thcontrex}
Существует потенциал $q$, принимающий чисто мнимые значения $q(x)\in i\RR$ при $x\in\RR_+$, такой, что $|q|\to+\infty$ при $x\to+\infty$,
но оператор $L_0$ не имеет расширений с компактной резольвентой.
\end{Theorem}

Для этого потенциала $\beta-\alpha=\pi$, и, очевидно, он удовлетворяет условию Молчанова. Оператор $L_U$ с краевыми условиями Дирихле
$U(y)=y(0)$ имеет ограниченную
резольвенту, по меньшей мере в левой полуплоскости \cite[Л.2]{Lidsky}, \cite{Birger}, но она не является вполне непрерывным оператором.

Следующая теорема даёт достаточное условие компактности резольвент операторов с потенциалами, удовлетворяющими $\RR^-$--условию с $\beta-\alpha>\pi$.
При этом теряется
свойство секториальности самих рассматриваемых операторов, в частности, числовой образ $L_U$ может заметать всю комплексную
плоскость \cite{Iskin}.
\begin{Theorem}
\label{thq32}
Пусть для некоторого $x_0>0$ при всех $x\ge x_0>0$ $|q(x)|\ge1$ и дополнительно:
\begin{itemize}
\item
$q\in AC_{loc}[x_0,+\infty)$,

\item для некоторого $0<\varkappa<\pi$
$$
-\pi+\varkappa<\arg q(x)<\pi-\varkappa,\quad x\ge x_0,
$$

\item
для некоторого $0<\delta<1$
$$
\left|\frac{q'(x)}{q^{3/2}(x)}\right|<4\delta\sin\frac{\varkappa}{2},\quad x\ge x_0.
$$
\end{itemize}
Тогда для компактности резольвенты $L_U$ достаточно, чтобы для любой бесконечной
системы непересекающихся интервалов $D_n\subset\RR_+$ равной длины
\begin{equation}
\label{eq32molchn}
\int\limits_{D_n}|q(x)|^{1/2}\,dx\to+\infty\ \mbox{при}\ n\to\infty.
\end{equation}
\end{Theorem}

Дальнейшее изложение работы: раздел \ref{secprooves} мы посвятим доказательствам теорем
\ref{thcriterion}, \ref{thcontrex} и \ref{thq32},
а в разделе \ref{sectnecessity} вернёмся к теореме \ref{thnecessn2}, которую сформулируем и докажем
для общих дифференциальных операторов порядков $n\ge 2$.

\section{Доказательство теорем \ref{thcriterion}, \ref{thcontrex} и \ref{thq32}}
\label{secprooves}
Так как ни условие Молчанова, ни факт абсолютной непрерывности резольвенты не зависят ни от сдвигов потенциала на постоянную, ни от значений потенциала на
конечном промежутке $[0,x_0]$, не ограничивая общности, $q_0=0$, $x_0=0$.

Доказательство критерия Молчанова \cite{Molchanov} основано на критерии Реллиха компактности резольвенты самосопряжённого положительно определённого
оператора \cite[Гл.II,\S24,Теор.11]{NaimarkBk}, который мы обобщим на случай m-секториальных операторов. Соответствующие определения см. \cite[Гл.V,\S3,п.10]{Kato}:
\begin{Lemma}
\label{lmrellich}
Пусть A --- m-секториальный оператор в гильбертовом пространстве $\mathfrak{H}$ с областью определения $\mathscr{D}_A\subset\mathfrak{H}$.

Оператор $A$ имеет компактную резольвенту тогда и только тогда, когда множество всех векторов $\varphi\in\mathscr{D}_A$,
удовлетворяющих условию
$$
\re(A\varphi,\varphi)\le1,
$$
компактно.
\end{Lemma}
{\noindent\bf Доказательство.} С оператором $A$ свяжем плотно определённую замкнутую секториальную полуторалинейную
форму $\mathfrak{a}$ \cite[Гл.VI,\S2,Теор.2.7]{Kato}, порождающим ядром которой будет $\mathscr{D}_A$.
По $\mathfrak{a}$ построим симметрическую форму $\mathfrak{t}=\re\mathfrak{a}$: для $u,v\in\mathscr{D}_\mathfrak{a}$ (из области определения $\mathfrak{a}$),
положим
$$
\mathfrak{t}[u,v]=\frac{1}{2}(\mathfrak{a}[u,v]+\overline{\mathfrak{a}[v,u]}).
$$

Форма $\mathfrak{t}$ замкнута на области $\mathscr{D}_\mathfrak{t}=\mathscr{D}_\mathfrak{a}$, плотно определена, неотрицательна;
по второй теореме о представлении \cite[Гл.VI,\S2,Теор.2.23]{Kato} с ней ассоциирован самосопряжённый оператор $T\ge0$,
$\mathscr{D}_\mathfrak{t}=\mathscr{D}(T^{1/2})\supset\mathscr{D}_T$, где область $\mathscr{D}_T$ определения $T$ является порождающим ядром $\mathfrak{t}$, и
$$
\mathfrak{t}[u,v]=(T^{1/2}u,T^{1/2}v),\quad u,v\in\mathscr{D}_\mathfrak{t}.
$$
Область $\mathscr{D}_A$ также является ядром $\mathfrak{t}$. Ввиду того, что область $\mathscr{D}_A$ --- ядро $\mathfrak{a}$,
для любого $u\in\mathscr{D}_\mathfrak{t}=\mathscr{D}_\mathfrak{a}$ найдём последовательность $u_n\in\mathscr{D}_A$ такую, что
$u_n\to u$ и $\mathfrak{a}[u_n,u_n]\to\mathfrak{a}[u,u]$ при $n\to\infty$. Очевидно, $\mathfrak{t}[u_n,u_n]=\re\mathfrak{a}[u_n,u_n]$ сходится,
а ввиду замкнутости формы $\mathfrak{t}$, предел $\lim\mathfrak{t}[u_n,u_n]=\mathfrak{t}[u,u]$.

Резольвенты обоих операторов $A$ и $T$ компактны или не компактны одновременно \cite[Гл.VI,\S3,Теор.3.3]{Kato}.

Последующие рассуждения используют критерий Реллиха для $T$.

Докажем необходимость. Пусть $\Phi=\set{\varphi\in\mathscr{D}_A}{\re(A\varphi,\varphi)\le1}$. Так как $\mathscr{D}_T$ --- ядро $\mathfrak{t}$,
для любого $\varepsilon>0$ и любого $\varphi\in\Phi$ найдём $u\in\mathscr{D}_T$ такое, чтобы
$|\mathfrak{t}[u,u]-\mathfrak{t}[\varphi,\varphi]|<\varepsilon$ и $|u-\varphi|<\varepsilon$. Из компактности резольвенты
$A$ следует компактность резольвенты $T$, а так как
$(Tu,u)=\mathfrak{t}[u,u]\le1+\varepsilon$, делаем вывод о компактности всей совокупности $\{u\}$, а ввиду произвольности $\varepsilon$
--- вывод о компактности $\Phi$.

Докажем достаточность. Пусть $U=\set{u\in\mathscr{D}_T}{(Tu,u)\le1}$.
Так как $\mathscr{D}_A$ --- ядро $\mathfrak{t}$,
для любого $\varepsilon>0$ и любого $u\in U$ найдём $\varphi\in\mathscr{D}_A$, чтобы
$|\mathfrak{t}[u,u]-\mathfrak{t}[\varphi,\varphi]|<\varepsilon$ и $|u-\varphi|<\varepsilon$.
Так как $\re(A\varphi,\varphi)=\mathfrak{t}[\varphi,\varphi]\le1+\varepsilon$, делаем вывод о компактности всей системы
$\{\varphi\}$, а ввиду произвольности $\varepsilon>0$ --- о компактности $U$, следовательно $T$ имеет компактную резольвенту,
а следовательно, и $A$.\qquad$\Box$
\bigskip

{\noindent\bf Доказательство теоремы \ref{thcriterion}.} Необходимость следует из теоремы \ref{thnecessn2}, остановимся только на достаточности.

Не ограничивая общности, пусть $-\pi<\alpha\le0\le\beta<\pi$ (в противном случае расширим границы сектора изменения потенциала, оставаясь в условиях теоремы).

Покажем, что у $L_U$ существует резольвента в некоторой области.

В условиях теоремы для однородного уравнения $l(y)=\lambda y$ реализуется случай предельной точки \cite[Сл.Теор.4]{Birger},
и в некотором секторе $\Lambda\subset\CC$ операторы $L_{D}$ и $L_N$ с краевыми условиями Дирихле $U(y)=y(0)$ и Неймана $U(y)=y'(0)$
имеют ограниченные резольвенты $R_{D,\lambda}$ и $R_{N,\lambda}$, соответственно \cite[Th.4.1]{Brown}.

Из явного выражения соответствующих резольвент через функции Грина следует, что у $L_U-\lambda$
при $\lambda\in\Lambda$ существует по меньшей мере правый обратный оператор
$$
R_\lambda=w_D(\lambda)R_{D,\lambda}+w_N(\lambda)R_{N,\lambda},\quad w_D(\lambda)+w_N(\lambda)=1,
$$
где $w_D$, $w_N$ --- мероморфные в $\Lambda$ функции. Т.е. по меньшей мере в некоторой подобласти $\lambda\in\Lambda_0\subset\Lambda$, не содержащей
полюсов $w_D$, $w_N$, оператор $L_U-\lambda$
сюръективен. Ввиду реализации случая предельной точки, сопряжённое дифференциальное выражение и сопряжённые краевые условия
задают сопряжённый оператор $(L_U-\lambda)^*$, который аналогично оказывается сюръективным, а значит $L_U-\lambda$
взаимно однозначно отображает $\mathscr{D}_U$ на все $L_2(\RR_+)$ при $\lambda\in\Lambda_0$, и резольвента определена однозначно.

При каждом $\lambda\in\Lambda$ дефект оператора $L_0-\lambda$ равен 1 (образ $L_0-\lambda$ ортогонален единственному решению из $L_2(\RR_+$) уравнения
$y''=(\overline{q(x)}-\overline{\lambda})y$), резольвенты любых двух расширений
$L_0$ в каждой фиксированной точке $\lambda\in\Lambda$, где они одновременно существуют, отличаются не более чем на одномерный оператор. Следовательно, они компактны
или не компактны одновременно. И дальнейшее доказательство достаточно провести лишь для оператора $L_D$ с краевыми условиями Дирихле.

Из существования резольвенты для $\lambda\in\Lambda$ следует замкнутость $L_D$.

Обозначим область определения $L_D$ через $\mathscr{D}_D$, также рассмотрим линейное многообразие $\mathscr{D}_{D0}$, включающее элементы $\mathscr{D}_D$
с компактным носителем:
\begin{equation}
\label{eqDDdefinit}
\mathscr{D}_D=\set{y\in\mathscr{D}}{y(0)=0},\quad
\mathscr{D}_{D0}=\set{y\in\mathscr{D}_D}{\exists x_0>0\ \forall x\ge x_0\ y(x)=0}.
\end{equation}
Область $\mathscr{D}_{D0}$ является порождающим ядром $L_D$ в том смысле, что замыкание ограничения $L_D$ на $\mathscr{D}_{D0}$ совпадает
с $L_D$ \cite[Л.6]{Lidsky}.

Положим $\theta=-(\alpha+\beta)/2$. Оператор $M=e^{i\theta}L_D$ m-секториальный.
Ввиду существования резольвенты он максимальный замкнутый, а также для любого $y\in\mathscr{D}_{D0}$
$$
(My,y)=e^{i\theta}\int\limits_0^{+\infty}|y'(x)|^2\,dx+\int\limits_0^{+\infty}e^{i\theta}q(x)|y(x)|^2\,dx,\quad |\arg(My,y)|\le\frac{\beta-\alpha}{2}<\frac{\pi}{2}.
$$

Далее рассуждения в существенном повторяют рассуждения Молчанова \cite{Molchanov}, тем не менее приведём их здесь, они пригодятся
при доказательстве теоремы \ref{thq32}.

Допустим, условие Молчанова выполнено, но резольвента $M$ не компактна, возьмём некомпактную последовательность
$Y=\sint{y_n}{n=1}{\infty}\subset\mathscr{D}_{D}$,
чтобы $\re(My_n,y_n)\le1$, $n\in\NN$ ввиду Леммы \ref{lmrellich}. Так как $\mathscr{D}_{D0}$ --- порождающее ядро $M$, не ограничивая общности,
все $y_n\in\mathscr{D}_{D0}$. Действительно, достаточно для любого $\varepsilon>0$ и любого $y_n\in Y$ найти $\hat y_n\in\mathscr{D}_{D0}$, чтобы
$\|y_n-\hat y_n\|<\varepsilon$ и $|\re(M\hat y_n,\hat y_n)-\re(My_n,y_n)|<\varepsilon$. Система
$\sint{\hat y_n}{n=1}{\infty}$ не может быть компактной при всех $\varepsilon>0$, найдётся $\varepsilon_0>0$, при котором $\sint{\hat y_n}{n=1}{\infty}$
не будет компактной. При этом $\re(M\hat y_n,\hat y_n)\le1+\varepsilon_0$.
Подходящая нам система будет состоять из нормированных элементов $\sint{\hat y_n/\sqrt{1+\varepsilon_0}}{n=1}{\infty}$.

Найдётся $\varepsilon_0>0$ такое, что для любого $T>0$ существует $y_T\in Y$ такое, что
\begin{equation}
\label{eqyTeps0}
\int\limits_T^{+\infty} |y_T(x)|^2\,dx\ge\varepsilon_0.
\end{equation}
Это следует из компактности срезок $Y_T=\sint{y_{n,T}}{n=1}{\infty}$, где $y_{n,T}(x)=y_n(x)$ при $x\in[0,T]$, $y_{n,T}(x)=0$ при $x>T$ \cite[Л.8]{Lidsky}. Если бы
для всякого $\varepsilon>0$ нашлось $T_0=T_0(\varepsilon)$ и для всякого $n\in\NN$ интегралы
$$
\int\limits_{T_0}^{+\infty} |y_n(x)|^2\,dx<\varepsilon,
$$
то и сама последовательность $Y$ оказалась компактной (см. заключительные рассуждения доказательства \cite[Теор.4]{Lidsky}).

Возьмём $d=(\varepsilon_0^{1/2}\cos^{1/2}\theta)/4>0$, произвольное $T>0$, соответствующее $y_T$, разобьём луч $[T,+\infty)$ отрезками $D_n$, $n\in\NN$ равной длины $d$,
хотя бы на одном из них,
который мы обозначим $D_T$,
\begin{gather}
\label{eqlgdtinq}
\re\left(
e^{i\theta}\int\limits_{D_T}|y_T'(x)|^2\,dx+\int\limits_{D_T}e^{i\theta}q(x)|y_T(x)|^2\,dx
\right)\le\frac{1}{\varepsilon_0}\int\limits_{D_T} |y_T(x)|^2\,dx,\\
\int\limits_{D_T} |y_T(x)|^2\,dx>0.
\label{eqlgdtinq2}
\end{gather}
Если бы неравенства противоположные \eqref{eqlgdtinq} были справедливы для всех $D_n$, $n\in\NN$, на которых выполнено \eqref{eqlgdtinq2},
суммируя их, и применяя \eqref{eqyTeps0}, получили бы противоречащую
оценку $\re(My_T,y_T)>1$.

Проведём нормировку $y_T$, полагая
\begin{equation}
\label{eqvtdefnts}
v_T=y_T \frac{d^{1/2}}{\Bigl(\int\limits_{D_T} |y_T(x)|^2\,dx\Bigr)^{1/2}},
\end{equation}
тогда, учитывая, что длина $|D_T|=d$,
\begin{equation}
\label{eqvTmain}
\re\left(
e^{i\theta}\int\limits_{D_T}|v_T'(x)|^2\,dx+\int\limits_{D_T}e^{i\theta}q(x)|v_T(x)|^2\,dx
\right)\le\frac{d}{\varepsilon_0},\quad\frac{1}{|D_T|}
\int\limits_{D_T}|v_T(x)|^2\,dx=1.
\end{equation}

В частности,
$$
\re\left(
e^{i\theta}\int\limits_{D_T}|v_T'(x)|^2\,dx
\right)\le\frac{d}{\varepsilon_0},\mbox{ т.е. }
\int\limits_{D_T}|v_T'(x)|^2\,dx\le\frac{d}{\varepsilon_0\cos\theta}.
$$

Для любых двух $x_1,x_2\in D_T$,
\begin{gather*}
\bigl||v_T(x_1)|-|v_T(x_2)|\bigr|\le
\bigl|v_T(x_1)-v_T(x_2)\bigr|\le\\
\le\int\limits_{D_T}|v_T'(x)|\,dx\le d^{1/2}\Bigl(
\int\limits_{D_T}|v_T'(x)|^2\,dx\Bigr)^{1/2}\le \frac{1}{4},
\end{gather*}
так как $|v_T|$ непрерывна на $D_T$, а интегральное среднее от $|v_T|^2$  равно 1 \eqref{eqvTmain},
найдётся $x_0\in D_T$, что $|v_T(x_0)|=1$. Из
полученной оценки делаем вывод, что $|v_T(x)|\ge3/4$ для всех $x\in D_T$.

Снова обращаемся к \eqref{eqvTmain}, откуда:
$$
\int\limits_{D_T}\re(e^{i\theta}q(x))|v_T(x)|^2\,dx
\le\frac{d}{\varepsilon_0},
$$
но $|\arg(e^{i\theta}q(x))|\le(\beta-\alpha)/2$, так что
$$
\re(e^{i\theta}q(x))\ge|q|\cos\frac{\beta-\alpha}{2},
$$
стало быть,
$$
\int\limits_{D_T}|q(x)|\,dx\le\frac{d}{\varepsilon_0}\frac{4^2}{3^2}\sec\frac{\beta-\alpha}{2},
$$
что, ввиду произвольности $T>0$, противоречит условию Молчанова для $q$. Полученное противоречие завершает
доказательство.\qquad$\Box$
\bigskip

{\noindent\bf Доказательство теоремы \ref{thcontrex}.} Построим контрпример, который и будет служить доказательством.

Разобьём полуось на полусегменты длины $\pi$, полагая $I_k=[\pi(k-1);\pi k)$, $k\in\NN$, на каждом
определяя $q_k(x)=ikr_{n_k}(x-\pi(k-1))$, $x\in I_k$, где
$$
r_n(t)=(-1)^{j+1},\ t\in\Bigl[\pi\frac{j-1}{n},\pi\frac{j}{n}\Bigr),\ j=1,\ldots,n;\quad t\in[0,\pi).
$$
Процедура выбора $\sint{n_k}{k=1}{\infty}$ будет ясна позже.

Полагая $q\vint{I_k}{}=q_k$, определим $q$ для всех $x\in\RR_+$. Очевидно, $|q(x)|\to+\infty$, когда
$x\to+\infty$.

В качестве краевых условий возьмём условия Дирихле $U(y)=y(0)$. Для упрощения записи сам оператор
обозначим через $L$.

Для однородного уравнения $y''(x)=(q(x)-\lambda)y(x)$ реализуется случай предельной точки, и, как уже отмечалось,
оператор $L$ имеет ограниченную резольвенту в левой полуплоскости \cite[Л.2]{Lidsky}, \cite{Birger}.
Покажем, как за счёт выбора $\sint{n_k}{k=1}{\infty}$ можно добиться того, что резольвента не будет компактным оператором.

Мы покажем, что для каждого $k\in N$ существует такое $n_k\in N$, и такая $y_k\in\mathscr{D}_{U}$,
равная нулю всюду вне $I_k$, что $\|y_k\|=1$ и $\|Ly_k\|<2^{13}$. Откуда и будет следовать некомпактность резольвенты
оператора $L$ с соответствующим набором $\sint{n_k}{k=1}{\infty}$.

Эта задача локальна для каждого $I_k$, в связи с чем достаточно для любого $u>0$ найти такое $n\in\NN$, чтобы для
оператора
$$
L_n(u)y=-y''+iur_ny
$$
в $L_2[0,\pi]$, заданного на области
$$
\mathscr{E}=\set{y\in L_2[0,\pi]}{y,y'\in AC[0,\pi],\ y''\in L_2[0,\pi],\ y(0)=y(\pi)=0},
$$
существовало $y\in\mathscr{E}$,
$y(0)=y'(0)=y(\pi)=y'(\pi)=0$, $\|y\|=1$ и $\|L_n(u)y\|<2^{13}$.
Здесь и далее фигурируют нормы пространства $L_2[0,\pi]$.

Вместе с операторами $L_n(u)$ введём в рассмотрение оператор $L_0$ с той же областью определения $\mathscr{E}$.
$$
L_0y=-y'',
$$
через $\mu_j=j^2$, $j\in\NN$ обозначим собственные значения $L_0$.

Оператор $B_ny=r_ny$ ограничен в $L_2[0,\pi]$, $\|B_n\|=1$, и
$B_n\stackrel{w}{\rightarrow}0$ при $n\to\infty$ в смысле слабой операторной сходимости,
что легко проверяется на индикаторах $\chi_{[a,b]}$ отрезков $[a,b]\subset [0,\pi]$,
линейные комбинации которых плотны в $L_2[0,\pi]$.

Покажем, что при фиксированном $u\ge0$ и $n\to\infty$ собственные значения $L_n(u)$ поточечно сходятся к $\sint{\mu_j}{j=1}{\infty}$:
для любого $\mu_j$ существует последовательность
собственных значений операторов $L_n(u)$: $\lambda_{j,n}\to \mu_j$, и для любого компакта $\mathcal{C}$, не содержащего точек
$\sint{\mu_j}{j=1}{\infty}$, найдётся
$n_0>0$, такое, что при всех $n>n_0$ компакт $\mathcal{C}$ не будет содержать собственных значений операторов $L_n(u)$.

Доказательство будем проводить индукцией по $u$.

При $u=0$ утверждение верно. Покажем, что если оно верно для некоторого $u_0\ge0$, то остаётся верным и для всех $u\in[u_0,u_0+1]$.

Пусть $u\in[u_0,u_0+1]$. Возьмём целое $j_0>\max\{u_0+1/2,\ u_0/2+13/8\}$.

Через $\Gamma_j$, $j\in\NN$ обозначим замкнутый круг с центром $\mu_j$. При $j=1,\ldots,j_0-1$ пусть радиусы кругов равны $5/4$, а при $j\ge j_0$ --- равны $u$.

Собственные значения $L_n(u)=L+iuB_n$ лежат в объединении замкнутых кругов с центрами в $\mu_j$ радиуса $u$.
При $j\ge j_0$ это непересекающиеся между собой круги $\Gamma_j$.

Так как $L_n(u)=L_n(u_0)+i(u-u_0)B_n$, из предположения индукции следует существование $n_0>0$ такого, что при всех $n>n_0$ и $j=1,\ldots,j_0-1$
каждый замкнутый круг с центром $\mu_j$ и радиусом $1/4$ содержит хотя бы одно собственное значение $L_n(u_0)$. Так как $0\le u-u_0\le1$, то для $j=1,\ldots,j_0-1$
замкнутые круги $\Gamma_j$ радиуса $5/4$ будут содержать хотя бы по одному собственному значению $L_n(u)$ и не будут пересекаться ни между собой,
ни с кругами $\Gamma_j$, $j\ge j_0$.

Мы нашли для каждого $j\in\NN$ изолированный круг $\Gamma_j$ с центром в $\mu_j$, содержащий по меньшей мере одно собственное значение $L_n(u)$ при всех $n>n_0$.
Возьмём произвольный
такой круг $\Gamma_j$ и рассмотрим в нем последовательность собственных значений $\lambda_n$ операторов $L_n(u)$. Пусть $\lambda_0\in\Gamma_j$ --- некоторая предельная точка,
покажем, что она совпадает с $\mu_j$. Предположим противное.

Пусть при $l\to\infty$ подпоследовательность $\lambda_{n_l}\to\lambda_0$, а $f_{n_l}$ --- нормированные собственные функции $L_{n_l}(u)$, $\|f_{n_l}\|=1$. Тогда
$$
(L_0-\lambda_{n_l})f_{n_l}=-iuB_{n_l}f_{n_l},\quad \|-iuB_{n_l}f_{n_l}\|=u.
$$

Положим $h_{n_l}=(L_0-\lambda_0)f_{n_l}$, $\|h_{n_l}\|=u$. Так как $\lambda_0\ne\mu_j$, тогда оператор $(L_0-\lambda_0)^{-1}$ существует и компактен,
компактной будет и последовательность $f_{n_l}$, которую, не ограничивая общности, можно считать сходящейся $f_{n_l}\to f$, $\|f\|=1$, следовательно,
$B_{n_l}f_{n_l}\stackrel{w}{\rightarrow}0$.

Положим $g=(L_0-\overline{\lambda_0})^{-1}f$, тогда
\begin{gather*}
\|f\|^2=\lim_{l\to\infty}(f_{n_l},(L_0-\overline{\lambda_0})g)=\lim_{l\to\infty}((L_0-\lambda_0)f_{n_l},g)=\\
=\lim_{l\to\infty}\bigl\{((L_0-\lambda_{n_l})f_{n_l},g)+iu(B_{n_l}f_{n_l},g)\bigr\}=
\lim_{l\to\infty}((L_{n_l}(u)-\lambda_{n_l})f_{n_l},g)=0.
\end{gather*}
Мы воспользовались тем, что $(\lambda_{n_l}-\lambda_0)f_{n_l}\to0$ и $(B_{n_l}f_{n_l},g)\to0$.

Но $\|f\|=1$, мы пришли к противоречию, следовательно, $\lambda_n\to\mu_j$. Аналогично доказывается, что любой компакт, не содержащий
собственных значений $L_0$, при всех достаточно больших $n>n_0$ не будет содержать собственных значений $L_n(u)$.

Положим $\mathscr{E}_0=\set{y\in \mathscr{E}}{y'(0)=y'(\pi)=0}$.

Покажем, что для любого $u>0$ можно подобрать $n_0>0$ и $y_{n_0}\in \mathscr{E}_0$, $\|y_{n_0}\|=1$ такие, что
$\|L_{n_0}(u)y_{n_0}\|<2^{13}$, что завершит построение контрпримера.

Обозначим $H_n=L_n(u)\mathscr{E}_0=L_2[0,\pi]\ominus <y_{n,1},y_{n,2}>$, где $y_{n,1}$ и $y_{n,2}$ --- ФСР
однородного уравнения
$$
-y''-iur_ny=0.
$$
Здесь $<y_{n,1},y_{n,2}>$ --- двумерная плоскость, натянутая на $y_{n,1}$ и $y_{n,2}$.

Через $R_n(u)$ обозначим обратный оператор к $L_n(u)$. Его существование и оценка $\|R_n(u)\|\le 1$ следуют из несложной выкладки:
$$
\|My\|\|y\|\ge\bigl|(My,y)
\bigr|=|(Ly,y)\pm i(uB_ny,y)
\bigr|\ge(Ly,y)\ge \|y\|^2,
$$
где $M=L_n(u)$, либо $M=L_n(u)^*$

Сужение $R_n(u)$ на $H_n$ обозначим через $\overset{\circ}{R}_n(u)$, оценим его норму как оператора из $H_n$ в $L_2[0,\pi]$

$$
\|\overset{\circ}{R}_n(u)\|^2=\sup_{f\perp<y_{n,1},y_{n,2}>}\frac{\|R_n(u)f\|^2}{\|f\|^2}\ge \min_{\mathscr{L},\,\dim\mathscr{L}=2}\ \max_{f\perp\mathscr{L}}
\frac{\bigl((R_n(u))^*R_n(u)f,f\bigr)}{\|f\|^2}=s_3^2,
$$
минимум берётся по всевозможным двумерным подпространствам $L_2[0,\pi]$, величина $s_3$ --- 3-е сингулярное число компактного оператора $R_n(u)$. Последнее
равенство следует из \cite[Гл.II,\S1]{GhbKrn}.

Принимая во внимание $1\ge \|R_n(u)\|=s_1\ge s_2$, c помощью леммы Вейля оценим
$$
s_3s_2s_1\ge|\lambda_1||\lambda_2||\lambda_3|,\ \mbox{откуда}\ s_3\ge|\lambda_1||\lambda_2||\lambda_3|,
$$
где $\lambda_j$, $j=1,2,3$ --- три максимальных по модулю собственных значения $R_n(u)$ --- обратные величины трёх минимальных по модулю собственных
значений $L_n(u)$, но последние при $n\to\infty$
сходятся к $\mu_j$, следовательно, при некотором большом $n_0$ каждое $|\lambda_j|>1/\mu_4=1/2^4$. Стало быть, для $R_{n_0}(u)$ величина $s_3>1/2^{12}$, и
$$
1\ge\|\overset{\circ}{R}_{n_0}(u)\|>1/2^{12},
$$
но тогда найдётся $f_{n_0}\in H_n$, $f_{n_0}\ne0$ такое, что $\|R_{n_0}(u)f_{n_0}\|\ge 1/2^{13}\|f_{n_0}\|$.

Полагая $y_{n_0}=R_{n_0}(u)f_{n_0}/\|R_{n_0}(u)f_{n_0}\|$, найдём искомый элемент $\mathscr{E}_0$.\qquad$\Box$
\bigskip

{\noindent\bf Доказательство теоремы \ref{thq32}.} Извлечём корень из $q$, полагая $p(x)=\sqrt{q(x)}$ при $x\ge0$,
выбрав ветвь так, чтобы $|\arg p(x)|\le \pi/2-\varkappa/2$. Положим, следуя \cite{T023Comp},
$$
\rho(x)=\re p(x)-\frac{1}{2}\left|
\frac{p'(x)}{p(x)}
\right|,\quad x\ge0.
$$

Проведём оценку, положив $C_0=(1-\delta)\sin\varkappa/2>0$:
\begin{equation}
\label{eqrhoestim}
\rho(x)=\re p(x)-\frac{|p(x)|}{4}\left|
\frac{q'(x)}{q^{3/2}(x)}
\right|\ge\re p(x)-|p(x)|\delta\sin\frac{\varkappa}{2}\ge C_0|p(x)|\ge C_0>0,
\end{equation}
так как $\re p(x)\ge |p(x)|\sin\varkappa/2$ ввиду $|\arg p(x)|\le \pi/2-\varkappa/2$, и $|p(x)|\ge1$.

Вследствие \cite[Теор.1]{T023Comp} у $L_U$ существует ограниченная резольвента в некоторой окрестности нуля $\lambda\in\Omega\subset\CC$.
Для $L_U$ реализуется определённый случай (аналог случая предельной точки).
Рассуждая как и при доказательстве теоремы \ref{thcriterion}, считаем, что $L_U$ задан краевыми
условиями Дирихле $U(y)=y(0)$ и определён на области $\mathscr{D}_D$ \eqref{eqDDdefinit}. Из \cite[Теор.1]{T023Comp} также заключаем,
что $\mathscr{D}_{D0}$ --- порождающее ядро для $L_U$.

Оператор $R_0=L_U^{-1}$ ограничен, и может быть продолжен как ограниченный оператор $\widetilde{R_0}:L_2(\RR_+,1/|q|)\to L_2(\RR_+)$ \cite[Теор.4]{T023Comp}.
Это эквивалентно ограниченности оператора $Rg=\widetilde{R_0}(pg)$ как обычного оператора в $L_2(\RR_+)$.
Так как из компактности $R$ вытекает
компактность $R_0$, наша цель --- доказать полную непрерывность оператора $R$.

Из ограниченности $R$ следует замыкаемость
$$
M=\frac{1}{p}L_U
$$
как оператора в $L_2(\RR_+)$, а так как $L_2(\RR_+)\subset L_2(\RR_+,1/|q|)$ плотно в метрике весового пространства,
замыкание $\overline M$ порождается ядром $\mathscr{D}_D$ и $\overline M^{-1}=R$.

Ввиду того, что
$\mathscr{D}_{D0}$ --- порождающее ядро для $L_U$, и $|p|>1$ делаем вывод, что $\mathscr{D}_{D0}$ --- порождающее ядро для $\overline M$.

Приступим непосредственно к доказательству компактности резольвенты $\overline M$.

Покажем, что $\overline M$ --- m-секториальный оператор. Для этого нам достаточно показать лишь секториальность формы $(My,y)$ для $y\in\mathscr{D}_{D0}$ (см.
аналогичные рассуждения доказательства теоремы \ref{thcriterion}).
\begin{gather}
\notag
(My,y)=\int\limits_0^{+\infty}y'(x)\Bigl(
-\frac{p'(x)}{p^2(x)}\overline{y(x)}+\frac{1}{p(x)}\overline{y'(x)}
\Bigr)+p(x)|y(x)|^2\,dx=\\
\label{eqMforminit}
=\int\limits_0^{+\infty}
p(x)|y(x)|^2+\overline{p(x)}\Bigl|\frac{y'(x)}{p(x)}\Bigr|^2
-\frac{p'}{p}\Bigl(\frac{y'(x)}{p(x)}\Bigr)\overline{y(x)}\,dx.
\end{gather}

Положим
$$
\xi(x)=-\frac{p'}{p}\Bigl(\frac{y'(x)}{p(x)}\Bigr)\overline{y(x)},\quad
|\xi(x)|\le\frac{1}{2}\Bigl|\frac{p'}{p}\Bigr|\Bigl(
\Bigl|\frac{y'(x)}{p(x)}\Bigr|^2+|y(x)|^2
\Bigr),
$$
для удобства представим в виде:
$$
\xi(x)=m(x)\frac{1}{2}\Bigl|\frac{p'}{p}\Bigr|\Bigl(
\Bigl|\frac{y'(x)}{p(x)}\Bigr|^2+|y(x)|^2
\Bigr)e^{i\varphi(x)},
$$
где $0\le m(x)\le1$, $\varphi(x)=\arg\xi(x)$. В итоге \eqref{eqMforminit} примет вид
\begin{gather}
\notag
(My,y)=\int\limits_0^{+\infty}
\Bigl(
p(x)+m(x)\frac{1}{2}\Bigl|\frac{p'}{p}\Bigr|e^{i\varphi(x)}
\Bigr)|y(x)|^2+\\
+\Bigl(
\overline{p(x)}+m(x)\frac{1}{2}\Bigl|\frac{p'}{p}\Bigr|e^{i\varphi(x)}
\Bigr)\Bigl|\frac{y'(x)}{p(x)}\Bigr|^2\,dx.
\label{eqMform2}
\end{gather}

Оценим
$$
\left|m(x)\frac{1}{2}\Bigl|\frac{p'}{p}\Bigr|e^{i\varphi(x)}
\right|\le\frac{1}{4}\Bigl|\frac{q'(x)}{q^{3/2}(x)}\Bigr||p(x)|\le
\delta\sin\frac{\varkappa}{2}|p(x)|.
$$

Обратимся к слагаемым в \eqref{eqMform2}. С учётом  $|\arg p(x)|\le\pi/2-\varkappa/2$ и полученной
оценки, из геометрических соображений следует, для некоторого $0<\epsilon\le\varkappa/2$:
$$
\left|\arg\Bigl(
p(x)+m(x)\frac{1}{2}\Bigl|\frac{p'}{p}\Bigr|e^{i\varphi(x)}
\Bigr)\right|\le\pi/2-\epsilon,\quad
\sin(\frac{\varkappa}{2}-\epsilon)\le\delta\sin\frac{\varkappa}{2},
$$
а значит, и $|\arg(My,y)|\le\pi/2-\epsilon$, и m-секториальность $\overline M$ доказана.

Далее нам пригодится оценка для для $y\in\mathscr{D}_{D0}$, которая немедленно выводится из \eqref{eqMform2} с учётом
\eqref{eqrhoestim}:
$$
\re (My,y)\ge\int\limits_0^{+\infty}C_0|p(x)|\Bigl(|y(x)|^2+\Bigl|\frac{y'(x)}{p(x)}\Bigr|^2\Bigr)\,dx.
$$

Последующие рассуждения в существенном повторяют доказательство теоремы \ref{thcriterion}. Остановимся лишь на
наиболее важных моментах.

Аналогично, предполагая, что резольвента $\overline M$ не компактна, найдём некомпактную последовательность
$Y=\sint{y_n}{n=1}{\infty}\subset\mathscr{D}_{D0}$, чтобы $\re (My_n,y_n)<1$ для всех $n\in\NN$.

Найдём $\varepsilon_0>0$ и для любого $T>0$ возьмём $y_T\in Y$, чтобы выполнялось \eqref{eqyTeps0}.

Для $d=C_0\varepsilon_0/4$, найдём отрезок $D_T$ длины $d$, чтобы на нем
\begin{gather*}
C_0\int\limits_{D_T}|p(x)||y_T(x)|^2+\frac{|y_T'(x)|^2}{|p(x)|}\,dx\le\frac{1}{\varepsilon_0}\int\limits_{D_T} |y_T(x)|^2\,dx,\\
\int\limits_{D_T} |y_T(x)|^2\,dx>0.
\end{gather*}

Проведём нормировку \eqref{eqvtdefnts}, в результате которой получим
\begin{equation}
\label{eqvT32main}
C_0\int\limits_{D_T}|p(x)||v_T(x)|^2+\frac{|v_T'(x)|^2}{|p(x)|}\,dx\le
\frac{d}{\varepsilon_0},\quad\frac{1}{|D_T|}
\int\limits_{D_T}|v_T(x)|^2\,dx=1.
\end{equation}

Для любых двух $x_1,x_2\in D_T$,
\begin{gather*}
\bigl||v_T(x_1)|^2-|v_T(x_2)|^2\bigr|\le
\bigl|v_T^2(x_1)-v_T^2(x_2)\bigr|\le\\
\le2\int\limits_{D_T}|v_T(x)||v_T'(x)|\,dx
=2\int\limits_{D_T}|p(x)|^{1/2}|v_T(x)|\frac{|v_T'(x)|}{|p(x)|^{1/2}}\,dx\le\\
\le\int\limits_{D_T}|p(x)||v_T(x)|^2+\frac{|v_T'(x)|^2}{|p(x)|}\,dx
\le \frac{d}{C_0\varepsilon_0}=\frac{1}{4},
\end{gather*}
с учётом \eqref{eqvT32main} делаем вывод, что $|v_T(x)|^2\ge3/4$ для всех $x\in D_T$.

Вновь обращаясь к \eqref{eqvT32main}, получаем оценку
$$
C_0\frac{3}{4}\int\limits_{D_T}|p(x)|\,dx\le C_0\int\limits_{D_T}|p(x)||v_T(x)|^2\,dx\le\frac{d}{\varepsilon_0},
$$
или финально
$$
\int\limits_{D_T}|q(x)|^{1/2}\,dx=\int\limits_{D_T}|p(x)|\,dx\le\frac{1}{3},
$$
что приводит к противоречию и завершает доказательство. \qquad$\Box$
\bigskip

Условия теоремы \ref{thq32} являются так же необходимыми для компактности резольвенты оператора $\overline M$,
что доказывается идентично \cite{Molchanov}. Это означает, что предложенный нами метод доказательства не позволяет
ослабить \eqref{eq32molchn}.

\section{Необходимое условие компактности резольвенты}
\label{sectnecessity}

Здесь мы докажем необходимое условие компактности резольвенты для обыкновенного дифференциального оператора произвольного
порядка $n\ge2$. Так как результаты предшествующих разделов нам не понадобятся, используем идентичные обозначения для схожих объектов.

Пусть заданы комплекснозначные функции $p_1\equiv C_0=\const$, $p_j\in L_{1,loc}(\RR_+)$, $j=2,\ldots,n$. Введём дифференциальное выражение
\begin{equation}
\label{eqdiffnexpr}
l(y)=\frac{d^n}{dx^n}y+\sum\limits_{j=1}^n p_j\frac{d^{n-j}}{dx^{n-j}}y
\end{equation}
и многообразие
\begin{gather*}
\mathscr{D}_0=\set{y\in L_2(\RR_+)}{y,y',\ldots,y^{(n-1)}\in AC_{loc}(\RR_+),\ l(y)\in L_2(\RR_+),\\
y(0)=\ldots=y^{(n-1)}(0)=0,\ \exists x_0>0\ \forall x\ge x_0\ y(x)=0}.
\end{gather*}

Рассмотрим дифференциальный оператор $L_0$ в $L_2(\RR_+)$ с областью определения $\mathscr{D}_0$, заданный дифференциальным выражением
\eqref{eqdiffnexpr}.

\begin{Theorem}
\label{thnecessity}
Для того, чтобы у оператора $L_0$ существовало расширение с компактной резольвентой, необходимо, чтобы для любой бесконечной
системы непересекающихся интервалов $D_k\subset\RR_+$ равной длины
$$
\int\limits_{D_k}\sum_{j=2}^n|p_j(x)|\,dx\to+\infty\ \mbox{при}\ k\to\infty.
$$
\end{Theorem}

Доказательству предпошлём одну Лемму.

\begin{Lemma}
\label{lminfdimntn}
Пусть на отрезке $[0,d]$ заданы непрерывные комплекснозначные функции $\phi_\nu$, $W_\nu$, $\nu=1,\ldots,n$, $n\in\NN$ такие, что
$W_n(0)=1$, и для любого $d_1>0$, $0<d_1\le d$, система функций $\sint{\phi_\nu}{\nu=1}{n}$ линейно независима на $[0,d_1]$.

Тогда оператор $R$ в $L_2[0,d]$, задаваемый выражением
$$
(Rg) (x)=\sum_{\nu=1}^n\phi_\nu(x)\int\limits_0^x W_\nu(\xi)g(\xi)\,d\xi,
$$
бесконечномерный (т.е. его образ не является конечномерным линеалом).
\end{Lemma}
{\noindent\bf Доказательство.} Построим последовательность $g_l\in L_2[0,d]$, $l\in\NN$, для которой
последовательность образов $f_l=Rg_l$ линейно независима.

Для этого сперва построим непересекающиеся интервалы
$$
I_l=(\alpha_l,\beta_l)\subset[0,d],\quad l=0,1,\ldots,\quad
0<\ldots<\alpha_1<\beta_1<\alpha_0<\beta_0,
$$
чтобы система $\sint{\phi_\nu}{\nu=1}{n}$ была линейно независимой
на каждом $I_l$. Достаточно найти один такой интервал $I_0=(\alpha_0,\beta_0)$, а дальнейшее построение провести индуктивно.

Выбираем произвольно $\beta_0>0$, $0<\beta_0<d$. Если $\sint{\phi_\nu}{\nu=1}{n}$ линейно независима на $(\beta_0/2,\beta_0)$, то построение завершено, если же нет, пусть
$\mathscr{L}\subset\CC^n$ --- максимальное линейное многообразие такое, что для любого $A=(A_\nu)\in\mathscr{L}$,
$\sum_\nu A_\nu\phi_\nu(x)\equiv0$ при всех $x\in(\beta_0/2,\beta_0)$. Ясно, что $1\le\dim\mathscr{L}\le n$.

Пусть $\alpha'$ --- точная нижняя грань таких $\alpha>0$,
что для интервала $(\alpha,\beta_0)$ найдётся $A=(A_\nu)\in\mathscr{L}$, чтобы $\sum_\nu A_\nu\phi_\nu(x)\equiv0$
при всех $x\in(\alpha,\beta_0)$.
Заведомо $0\le\alpha'\le\beta_0/2$.

Покажем, что $\alpha'>0$. Иначе, существует $A_k=(A_{k\nu})\in\mathscr{L}$, $\|A_k\|=1$, $k\in\NN$, и
$\sum_\nu A_{k\nu}\phi_\nu(x)\equiv0$ при всех $x\in(\beta_0/(k+1),\beta_0)$. Выделим сходящуюся подпоследовательность $A_{k_s}\to A_0=(A_{0\nu})\in\mathscr{L}$
при $s\to\infty$, $\|A_0\|=1$.

Для $s\ge s_0\ge1$ линейная комбинация $\sum_\nu A_{0\nu}\phi_\nu(x)=\sum_\nu(A_{0\nu}-A_{k_s\nu})\phi_\nu(x)$ при всех $x\in(\beta_0/(k_{s_0}+1),\beta_0)$,
устремляя $s\to\infty$, в ввиду произвольности $s_0$ и непрерывности всех $\phi_\nu$,
получим $\sum_\nu A_{0\nu}\phi_\nu(x)=0$ для всех $x\in[0,\beta_0]$, что противоречит независимости $\sint{\phi_\nu}{\nu=1}{n}$ на отрезке $[0,\beta_0]$.
Таким образом, $\alpha'>0$.

Положим $\alpha_0=\alpha'/2$. Тогда система $\sint{\phi_\nu}{\nu=1}{n}$ окажется линейно независимой на $I_0=(\alpha_0,\beta_0)$. В противном случае
для некоторого $B=(B_\nu)\in\CC^n$ линейная комбинация
$\sum_\nu B_\nu\phi_\nu(x)\equiv0$ при всех $x\in(\alpha_0,\beta_0)$, а значит, и для $x\in(\beta_0/2,\beta_0)$, т.е. $B\in\mathscr{L}$, что противоречит выбору
$\alpha'$.

Пусть $I_l=(\alpha_l,\beta_l)$, $l=0,1,\ldots$ --- система непересекающихся интервалов линейной независимости $\sint{\phi_\nu}{\nu=1}{n}$,
считаем, что $\beta_0$ выбрано таким образом, чтобы $W_n(x)\ne0$ при $x\in[0,\beta_0]$.

Для $l\in\NN$ построим $g_l\in L_2[0,d]$ следующим образом: вне $I_l$ она будет тождественна 0, а на самом $I_l$ руководствуемся следующей процедурой.

Система функций $\sint{W_\nu}{\nu=1}{n}$, будучи рассмотренной на $I_l$, может быть выражена через $1\le k\le n$ независимых $W_{s_j}$ (по крайней мере,
$W_n(x)\not\equiv0$ на $I_l$):
$$
W_\nu=\sum_{j=1}^k A_{\nu j}W_{s_j},
$$
выберем $g_l\in L_2(I_l)$, чтобы $\overline{g_l}$ не была ортогональна $W_{s_1}$, а если $k>1$ дополнительно чтобы $\overline{g_l}$
была ортогональна всем $\sint{W_{s_j}}{j=2}{n}$ в смысле метрики $L_2(I_l)$.

Положим $f_l=Rg_l$, $l\in\NN$. Тогда
$$
f_l(x)=
\begin{cases}
0 & x\in[0,\alpha_l],\\
\sum\limits_{\nu=1}^nA_{\nu 1}\phi_\nu(x)\int\limits_{\alpha_l}^{\beta_l} W_{s_1}(\xi)g_l(\xi)\,d\xi, & x\in[\beta_l,d].
\end{cases}
$$

Для нас существенно, что интеграл в формуле выше отличен от нуля, а линейная комбинация $\sum_\nu A_{\nu 1}\phi_\nu(x)\not\equiv0$ при $x\in I_j$ при всех
$0\le j<l$.

Таким образом, $f_l(x)\equiv0$ при $x\in I_j$, когда $j>l$, но $f_l(x)\not\equiv0$ при $x\in I_j$, когда $0\le j<l$.
Последовательность $f_l$ линейно независима. \qquad$\Box$
\bigskip

{\noindent\bf Доказательство теоремы \ref{thnecessity}.} Допустим, противное, существует расширение $L\supset L_0$ с компактной резольвентой, и для некоторой
системы $D_k=[a_k,b_k]\subset\RR_+$ непересекающихся интервалов равной длины $d=|D_k|$
\begin{equation}
\label{eqpjintCj}
\int\limits_{D_k}|p_j(x)|\,dx<C_j,\quad j=1,\ldots,n,\ k\in\NN.
\end{equation}

Так как $p_1\equiv \const$, нам удобно записать \eqref{eqpjintCj} для всех $p_j$. Не ограничивая общности, все $C_j>0$.

Уменьшим при необходимости длины интервалов, чтобы выполнялась оценка:
\begin{equation}
\label{eqeta}
\eta=\sum\limits_{j=1}^nC_jd^{j-1}<1.
\end{equation}

Мы придём к противоречию, если найдём последовательность $y_k\in\mathscr{D}_0$, равных нулю вне $D_k$,
$\|y_k\|=1$, $\|L_0y_k\|<C$ для любого $k\in\NN$.

Для этого рассмотрим операторы $L_{k,0}$ в $L_2[0,d]$, заданные дифференциальными выражениями $l_k$:
$l_k(y(x))=l(y(x+a_k))$, $x\in[0,d]$ на областях
\begin{gather*}
\widetilde{\mathscr{D}_{0k}}=\set{y\in L_2[0,d]}{y,y',\ldots,y^{(n-1)}\in AC[0,d],\ l_k(y)\in L_2[0,d],\\
y(0)=\ldots=y^{(n-1)}(0)=y(d)=\ldots=y^{(n-1)}(d)=0}.
\end{gather*}

Найдя последовательность $f_k\in\widetilde{\mathscr{D}_{0k}}$, $\|f_k\|=1$, $\|L_{k,0}f_k\|<C$ для произвольного $k\in\NN$,
мы решим задачу, взяв $y_k(x)=f_k(x-a_k)$ для $x\in D_k$, а вне $D_k$ продолжив $y_k$ нулём.

Рассмотрим специальные решения $\psi_{k,\nu}$, задачи Коши для каждого однородного уравнения
$l_k(y)=0$, $k\in\NN$, заданные начальными условиями
$$
\psi_{k,\nu}^{(j-1)}(0)=\delta_\nu^j,\quad j,\nu=1,\ldots,n,\ k\in\NN,
$$
где $\delta_\nu^j$ --- символ Кронекера.

Существование таких решений следует из \cite[Гл.V,\S 16,Теор.1]{NaimarkBk}. Все $\psi_{k,\nu}$ и их производные до $n-1$ порядка включительно,
абсолютно непрерывны на $[0,d]$.

Зафиксируем $k\in\NN$. Для любой $g\in L_2[0,d]$ уравнение $l_k(f)=g$ имеет единственное решение, заданное
условиями: $f_k^{(j)}(0)=0$,
$j=0,\ldots,n-1$  \cite[Гл.V,\S 16,Теор.1]{NaimarkBk}. Его можно представить в виде $f=R_kg$, где
$R_k$ --- оператор Коши --- ограниченный вольтерров оператор в $L_2[0,d]$:
\begin{equation}
\label{eqRkg}
f=(R_kg) (x)=
\int\limits_0^x\frac{1}{W(\xi)}
\begin{vmatrix}
\psi_{k,1}(\xi) & \psi_{k,2}(\xi) & \cdots & \psi_{k,n}(\xi)\\
\psi_{k,1}'(\xi) & \psi_{k,2}'(\xi) & \cdots & \psi_{k,n}'(\xi)\\
\vdots & \vdots & \ddots & \cdots \\
\psi_{k,1}^{(n-2)}(\xi) & \psi_{k,2}^{(n-2)}(\xi) & \cdots & \psi_{k,n}^{(n-2)}(\xi)\\
\psi_{k,1}(x) & \psi_{k,2}(x) & \cdots & \psi_{k,n}(x)
\end{vmatrix}
g(\xi)\,d\xi,
\end{equation}
где $W(\xi)$ --- вронскиан
$$
W(\xi)=
\begin{vmatrix}
\psi_{k,1}(\xi) & \psi_{k,2}(\xi) & \cdots & \psi_{k,n}(\xi)\\
\psi_{k,1}'(\xi) & \psi_{k,2}'(\xi) & \cdots & \psi_{k,n}'(\xi)\\
\vdots & \vdots & \ddots & \cdots \\
\psi_{k,1}^{(n-2)}(\xi) & \psi_{k,2}^{(n-2)}(\xi) & \cdots & \psi_{k,n}^{(n-2)}(\xi)\\
\psi_{k,1}^{(n-1)}(\xi) & \psi_{k,2}^{(n-1)}(\xi) & \cdots & \psi_{k,n}^{(n-1)}(\xi)\\
\end{vmatrix}
=\exp\Bigl(-\int\limits_0^\xi p_1(a_k+t)\,dt\Bigr)=e^{-\xi C_0}.
$$

Дифференцируя \eqref{eqRkg} $n-1$ раз, убеждаемся, что $f\in\widetilde{\mathscr{D}_{0k}}$ тогда и только тогда, когда:
$$
\int\limits_0^d e^{\xi C_0}
\begin{vmatrix}
\psi_{k,1}(\xi) & \psi_{k,2}(\xi) & \cdots & \psi_{k,n}(\xi)\\
\psi_{k,1}'(\xi) & \psi_{k,2}'(\xi) & \cdots & \psi_{k,n}'(\xi)\\
\vdots & \vdots & \ddots & \cdots \\
\psi_{k,1}^{(n-2)}(\xi) & \psi_{k,2}^{(n-2)}(\xi) & \cdots & \psi_{k,n}^{(n-2)}(\xi)\\
\psi_{k,1}^{(s)}(d) & \psi_{k,2}^{(s)}(d) & \cdots & \psi_{k,n}^{(s)}(d)
\end{vmatrix}
g(\xi)\,d\xi=0,\quad s=0,\ldots,n-1,
$$
другими словами, когда $g\in\overset{\circ}{\mathfrak{H}}_k=L_2[0,d]\ominus\mathfrak{K}_k$, где $\mathfrak{K}_k$ --- конечномерное подпространство $L_2[0,d]$,
натянутое на функции, комплексно сопряжённые к $W_{ks}(x)$:
$$
W_{ks}(x)= e^{x C_0}
\begin{vmatrix}
\psi_{k,1}(x) & \psi_{k,2}(x) & \cdots & \psi_{k,n}(x)\\
\psi_{k,1}'(x) & \psi_{k,2}'(x) & \cdots & \psi_{k,n}'(x)\\
\vdots & \vdots & \ddots & \cdots \\
\psi_{k,1}^{(n-2)}(x) & \psi_{k,2}^{(n-2)}(x) & \cdots & \psi_{k,n}^{(n-2)}(x)\\
\psi_{k,1}^{(s)}(d) & \psi_{k,2}^{(s)}(d) & \cdots & \psi_{k,n}^{(s)}(d)
\end{vmatrix}
,\quad s=0,\ldots,n-1,
$$
следовательно, $\dim\mathfrak{K}_k\le n$.

Оператор $L_{k,0}$ взаимно однозначно отображает $\widetilde{\mathscr{D}_{0k}}$ на $\overset{\circ}{\mathfrak{H}}_k$. Через
$\overset{\circ}{R}_k\subset R_k$ обозначим обратный, действующий взаимно однозначно из $\overset{\circ}{\mathfrak{H}}_k$ на $\widetilde{\mathscr{D}_{0k}}$.

Чтобы найти искомую последовательность $f_k$, достаточно показать, что нормы $\|\overset{\circ}{R}_k\|\ge C>0$ --- равномерно ограничены снизу для
всех $k\in\NN$. Тогда для каждого $k$ существует $u_k\in\overset{\circ}{\mathfrak{H}}_k$, $u_k\ne0$ такое, что для $v_k=\overset{\circ}{R}_ku_k$ верна
оценка $\|v_k\|\ge C/2\|u_k\|$, искомое $f_k=v_k/\|v_k\|$, $\|L_{k,0}f_k\|\le 2/C$.
$$
\|\overset{\circ}{R}_k\|^2=\sup_{g\in\overset{\circ}{\mathfrak{H}}_k}\frac{\|R_kg\|^2}{\|g\|^2}
\ge\min_{\substack{\mathfrak{L}\subset L_2[0,d],\\ \dim\mathfrak{L}\le n}}\ \max_{\substack{g\in L_2[0,d],\\ g\perp\mathfrak{L}}}
\frac{(R_k^*R_k^{\phantom{*}} g,g)}{(g,g)}=s_{n+1}^2(R_k).
$$
Минимум берётся по всевозможным линеалам $\mathfrak{L}\subset L_2[0,d]$ размерности не выше $n$.
Последнее равенство вытекает из \cite[Гл.II,\S1]{GhbKrn}, где $s_{n+1}(R_k)$ --- $n+1$-е сингулярное число оператора $R_k$.

Далее мы покажем, что из последовательности операторов $R_k$ можно выделить равномерно сходящуюся подпоследовательность $R_{k_m}\to R$, $m\to\infty$ к
бесконечномерному оператору $R$, следовательно, имеющему положительными все $s$--числа.
Сингулярные числа $s_l(R_{k_m})\to s_l(R)>0$, $l\in\NN$ \cite[Гл.II,\S2,Сл.2.3]{GhbKrn}.
В частности, при $m>m_0$, $\|\overset{\circ}{R}_{k_m}\|\ge s_{n+1}(R_{k_m})\ge s_{n+1}(R)/2=C>0$.

Рассматривая $D_k$ с номерами индексов, соответствующим $k_m$, не ограничивая общности, можем считать, что для всех $k\in\NN$,
$\|\overset{\circ}{R}_k\|\ge C>0$. Это завершит доказательство.

Покажем, как выделить равномерно сходящуюся подпоследовательность из $R_k$.

Проинтегрируем $s$ раз (для $s=1,\ldots,n$) тождество $l_k(\psi_{k,\nu})=0$:
\begin{equation}
\label{eqpsiknuviaints}
\psi_{k,\nu}^{(n-s)}(x)=\sum\limits_{r=0}^{s-1}\frac{x^r}{r!}\psi_{k,\nu}^{(n-s+r)}(0)-
\int\limits_0^x d\xi_1\int\limits_0^{\xi_1}\cdots\,d\xi_{s-1}\int\limits_0^{\xi_{s-1}}
\sum\limits_{j=1}^n \widetilde{p_{jk}}(\xi_s)\psi_{k,\nu}^{(n-j)}(\xi_s)\,d\xi_s,
\end{equation}
где $\widetilde{p_{jk}}(x)=p_j(x+a_k)$. При $s=1$ внутренние интегралы отсутствуют.

Обозначим $M_{k,\nu,l}=\max|\psi_{k,\nu}^{(l)}(x)|$, $x\in[0,d]$, $l=0,\ldots,n-1$. Используя \eqref{eqpsiknuviaints} и оценки \eqref{eqpjintCj}, заключаем:
$$
M_{k,\nu,n-s}\le e^d+d^{s-1}\sum\limits_{j=1}^nC_j M_{k,\nu,n-j},
$$
умножая обе стороны на $C_s$ и суммируя по $s=1,\ldots,n$, получим
$$
\sum\limits_{s=1}^{n}C_sM_{k,\nu,n-s}\le e^d\sum\limits_{s=1}^{n}C_s+\eta\sum\limits_{j=1}^nC_j M_{k,\nu,n-j},
$$
откуда, с учётом \eqref{eqeta},
$$
\sum\limits_{s=1}^{n}C_sM_{k,\nu,n-s}\le\frac{e^d}{1-\eta}\sum\limits_{s=1}^{n}C_s.
$$

Так как все $C_s>0$, эта оценка означает, что все величины $M_{k,\nu,l}$ равномерно ограничены по $k\in\NN$, $\nu=1,\ldots,n$, $l=0,\ldots,n-1$.

При фиксированных $\nu=1,\ldots,n$, $l=0,\ldots,n-2$ рассмотрим последовательность $\{\psi_{k,\nu}^{(l)}\}_{k\in\NN}$. Она предкомпактна в
равномерной метрике на $[0,d]$ по теорема Арцела--Асколи. Выберем последовательность номеров $k_m$, чтобы для всех $\nu=1,\ldots,n$, $l=0,\ldots,n-2$
подпоследовательности $\{\psi_{k_m,\nu}^{(l)}\}_{m\in\NN}$ были равномерно сходящимися на $[0,d]$.
$$
\psi_{k_m,\nu}^{(l)}\Rightarrow\phi_\nu^l,\quad m\to\infty,
$$
где $\phi_\nu^l$ --- непрерывные функции, $\nu=1,\ldots,n$, $l=0,\ldots,n-2$ --- индексы.

Для удобства обозначим $\phi_\nu=\phi_\nu^0$. Ввиду равномерной сходимости, из известной теоремы о предельном переходе
\cite[Гл.XVI,Теор.4]{Zorich} следует,
что функции $\phi_\nu$ дифференцируемы $n-2$--кратно, и $\phi_\nu^l=\phi_\nu^{(l)}$, $l=0,\ldots,n-2$.

Отмечая, что детерминант в определении $R_k$ зависит только от производных $\psi_{k,\nu}^{(l)}$ не выше порядка $n-2$, делаем вывод, что
$R_{k_m}\Rightarrow R$ при $m\to\infty$ в смысле равномерной операторной сходимости, где
\begin{equation}
\label{eqRdef}
(Rg) (x)=
\int\limits_0^xe^{\xi C_0}
\begin{vmatrix}
\phi_1(\xi) & \phi_2(\xi) & \cdots & \phi_n(\xi)\\
\phi_1'(\xi) & \phi_2'(\xi) & \cdots & \phi_n'(\xi)\\
\vdots & \vdots & \ddots & \cdots \\
\phi_1^{(n-2)}(\xi) & \phi_2^{(n-2)}(\xi) & \cdots & \phi_n^{(n-2)}(\xi)\\
\phi_1(x) & \phi_2(x) & \cdots & \phi_n(x)
\end{vmatrix}
g(\xi)\,d\xi,\quad g\in L_2[0,d].
\end{equation}

Оператор $R$ компактен как равномерный предел компактных операторов. Нам осталось доказать бесконечномерность $R$.

Отметим, что для любого
отрезка $[0,d_1]\subset[0,d]$ система функций $\sint{\phi_\nu}{\nu=1}{n}$ линейно независима на $[0,d_1]$.

Предполагая противное, существуют $A_\nu\in\CC$ такие, что $\sum_\nu A_\nu\phi_\nu(x)\equiv0$ для всех $x\in[0,d_1]$.
Обратимся к \eqref{eqpsiknuviaints}, выражая линейные комбинации $\sum_\nu A_\nu\psi_{k,\nu}^{(n-s)}$ для $s=n$.
Получим
\begin{gather*}
\sum\limits_{\nu=1}^n A_\nu\psi_{k,\nu}(x)=\sum\limits_{r=0}^{n-1}\frac{x^r}{r!}\sum\limits_{\nu=1}^nA_\nu\psi_{k,\nu}^{(r)}(0)-\\
-C_0\int\limits_0^x d\xi_1\int\limits_0^{\xi_1}\cdots\,d\xi_{n-1}\int\limits_0^{\xi_{n-1}}
\sum\limits_{\nu=1}^nA_\nu\psi_{k,\nu}^{(n-1)}(\xi_n)\,d\xi_n-\\
-\int\limits_0^x d\xi_1\int\limits_0^{\xi_1}\cdots\,d\xi_{n-1}\int\limits_0^{\xi_{n-1}}
\sum\limits_{j=2}^n \widetilde{p_{jk}}(\xi_n)\sum\limits_{\nu=1}^nA_\nu\psi_{k,\nu}^{(n-j)}(\xi_n)\,d\xi_n,
\end{gather*}
полагая $k=k_m$ и переходя к пределу при $m\to\infty$, получим равномерно по $x\in[0,d_1]$
$$
0\equiv\sum\limits_{r=0}^{n-1}\frac{x^r}{r!}\sum\limits_{\nu=1}^nA_\nu\delta_\nu^{r+1}(0)
+C_0
\sum\limits_{r=1}^{n-1}\frac{x^r}{r!}\sum\limits_{\nu=1}^nA_\nu\delta_\nu^r,
$$
или
$$
0\equiv A_1+\sum\limits_{r=1}^{n-1}\bigl(A_{r+1}+C_0A_r\bigr)\frac{x^r}{r!},
$$
откуда все $A_\nu=0$, $\nu=1,\ldots,n$.

Представим $R$ в виде
$$
(Rg) (x)=\int\limits_0^x\sum_{\nu=1}^n\phi_\nu(x)W_\nu(\xi)g(\xi)\,d\xi,\quad g\in L_2[0,d],
$$
где $W_\nu$ --- непрерывные функции --- алгебраические дополнения в разложении определителя в \eqref{eqRdef} по последней строке, помноженные на $e^{\xi C_0}$,
$W_n(0)=1$. Бесконечномерность $R$ вытекает из леммы \ref{lminfdimntn}.\qquad$\Box$
\bigskip

Автор выражает благодарность члену-корреспонденту РАН профессору Андрею Андреевичу Шкаликову за внимание, поддержку и живой интерес.

Работа выполнена при содействии РНФ, (грант 20-11-20261).

\end{document}